\documentclass[reqno]{amsart}
\usepackage{hyperref}
\usepackage{amsmath, amsthm}
\usepackage{graphics}
\usepackage{graphicx}

\begin{document}

\pagenumbering{arabic}
\pagestyle{plain}

\title{A chemotaxis-Navier-Stokes system
 with dynamical boundary conditions}

\author {Baili Chen}

\address{Baili Chen \newline
Department of Mathematics and Computer Science\\
Gustavus Adolphus College\\
Saint Peter, MN 56082, USA}
 \email{bchen@gustavus.edu}


\begin{abstract}
A chemotaxis-Navier-Stokes system is studied under dynamical boundary conditions in a bounded convex domain
$\Omega\in \mathbb{R}^3$ with smooth boundary.  This models the interaction of populations of swimming bacteria
with the surrounding fluid.  The existence of a global weak solution is proved using multiple layers of approximations
and Rothe's method for the time discretization.
 \end{abstract}

\maketitle

\numberwithin{equation}{section}
\newtheorem{theorem}{Theorem}[section]
\newtheorem{lemma}[theorem]{Lemma}
\newtheorem{remark}[theorem]{Remark}
\newtheorem{definition}[theorem]{Definition}
\allowdisplaybreaks

\section{Introduction}

The Chemotaxis-Navier-Stokes system

\begin{equation}  \label {e1.1}
\left\{
\begin{array}{lcl}
\partial_t c -\alpha \Delta c  + u \cdot\nabla c  + n f(c) =0   &\mbox{a.e.  in}  &\Omega\times(0,T), \\
\partial_t n  - \nabla\cdot  \left(\beta  \nabla n -g(n,\ c)  \nabla c \right)  + u  \cdot\nabla n =0 &\mbox{a.e.  in}  &\Omega\times(0,T),\\
\partial_t u  - \nabla\cdot  (  \xi  \nabla u)  +(u\cdot \nabla) u +\nabla p = n \nabla   \sigma  &\mbox{a.e.  in}  &\Omega\times(0,T), \\
\nabla \cdot u =0 &\mbox{a.e.  in}  &\Omega\times(0,T).
\end{array} \right.
\end{equation}

in a bounded convex domain with smooth boundary describes an oxygen-driven bacteria suspension swimming in an incompressible
fluid like water.

The above system, which was first introduced by Tuval et al.
\cite{tuv}, consists of three coupled equations: an equation for the concentration of oxygen 
$c$, an equation for the population density $n$ of the bacteria, and the Navier-Stokes equation describing the water flow $u$.

Equation  \eqref{e1.1}  has been studied by several authors (e.g. \cite{win1}, \cite{win2}, \cite{brau}).  In their works, the equations are endowed with Neumann or Robin boundary
conditions for both $c$ and $n$.

In this work, we extend the model by assuming dynamical boundary condition for oxygen concentration $c$, which is given by
\begin{equation}\label{e1.2}
\partial_t c =  \Delta_{\tau} c -b \partial_{\eta}c \quad   \mbox{on}\quad   \partial\Omega
\end{equation}

Here we assume there is an oxygen source acting on the boundary, which depends on the oxygen flux $\partial_{\eta}c$ across the
boundary.  The Laplace-Beltrami operator $\Delta_{\tau}$ on the boundary describes the diffusion of oxygen along the boundary. 
The derivation of dynamical boundary conditions is introduced in \cite{gold}.

To the best of our knowledge,  Chemotaxis-Navier-Stokes system with the above dynamical boundary condition has not been
addressed in the existing literature.

To complete the system, we introduce the following boundary conditions for bacterial and fluid field.
\begin{eqnarray}
\beta  \partial_{\eta}n =  g(n,\ c)  \partial_{\eta}c \quad   \mbox{on}\    \partial\Omega\label{e1.3}\\
u=0\quad   \mbox{on}\quad    \partial\Omega\label{e1.4}
\end{eqnarray}

together with the initial conditions
\begin{equation}\label{e1.5}
c(x,\ 0)=c_0(x),\ n(x,\ 0)=n_0(x),\ u(x,\ 0)=u_0(x)
\end{equation}

The aim of this paper is to derive the existence of a global weak solution of the system \eqref{e1.1} - \eqref{e1.5}, which describes the interaction between
the bacteria density $n$, the oxygen concentration $c$, the fluid velocity field $u$ and the associated pressure $p$ in a bounded 
convex domain $\Omega$ with smooth boundary.  We set the gradient of gravitational potential $\sigma$ to be constant (i.e. $\nabla\sigma \equiv $ const).
 
The existence of a global weak solution is proved by the discretization in time (Rothe's method).
This technique was also used in
previous papers (e.g. \cite{vaz},  \cite{vla},  \cite{nec},  \cite{kac}) to solve other types of problems.

To fix the notation, we denote by $H^m(\Omega)$ the standard Sobolev space in $L^2(\Omega)$ with derivative of order less than
or equal to $m$ in $L^2(\Omega)$.  Let $D(\Omega)$ be the space 
of $C^{\infty}$ function with compact support contained in $\Omega$.  The closure of  $D(\Omega)$ in $H^m(\Omega)$ is denoted by
$H_0^m(\Omega)$.

Let $\Upsilon$ be the space
\[  \Upsilon = \left\{ u\in D(\Omega),  \nabla\cdot u=0 \right\}  \]

The closure of $\Upsilon$ in $L^2(\Omega)$ and in $H_0^1(\Omega)$ are denoted by $H$ and $V$ respectively.

We denote by $L^r(0,\ T; \ X)$ the Banach space of all measurable functions 
\[v:\  [0,\ T] \rightarrow X\]
 with norm
\[
\begin{array}{ll}
 \   & \Vert v \Vert _{L^r(0,\ T; \ X)} = {\left( \int_0^T  \Vert v \Vert _{X}^r dt  \right)}^{\frac{1}{r}},  \ \mbox{for}\   1\le r< \infty  \\
\mbox{or}  & \Vert v \Vert _{  L^{\infty}  (0,\ T; \ X)  } = \mbox{ess  sup}_{0\le t\le T} \Vert v \Vert _{X},  \ \mbox{for}\   r=\infty.
\end{array}
 \]
 The trace of a function is denoted by the subscript $\tau.$  For example, $c_{\tau}$ denotes the trace of function $c$.
  
 Throughout this paper, we denote by $M$ and $C$ the constants whose values may be different at each occurrence.
 
 Before stating the main result, we make the following assumptions throughout this paper:
 \[
 \begin{array}{lll}
(H_1) &\  f(\cdot) \ in\  C^0(R)  & \mbox{with}  \  f_0\le  f(\cdot) \le f_1;\   f_0,\ f_1 \in R^{+}.\\
(H_2) &\  g(\cdot,\ \cdot) \ in\  C^1(R^2)  & \mbox{with}  \  | g(\cdot,\ \cdot) | \le g_1;\    g_1 \in R^{+}.\\
(H_3) &\   \alpha,\  \beta,  \  \xi,\  b \in R^{+}.  & \   
 \end{array}
 \]
 
 The main result of this paper is:
 \begin{theorem} \label{th1.1}
Suppose $(H_1)-(H_3)$ hold, $(c_0,\  n_0,\ u_0)\in (L^2(\Omega))^2\times H.$  Then there exists 
functions $(c,\ n)\in \left( L^{\infty}(0,\ T, \ L^2(\Omega))\bigcap L^2(0,\ T, \ H^1(\Omega))   \right)^2$ and $u\in 
 L^{\infty}(0,\ T, \  H)\bigcap L^2(0,\ T, \ V)$ such that $\left( c(0),\ n(0),\ u(0)    \right)=(c_0,\ n_0,\ u_0)$
 and
 
 \begin{equation} \label {e1.6}
\left\{
\begin{array}{ll}
\  &\int_{\Omega} \partial_t c(t)\phi_1 dx +  \alpha \int_{\Omega} \nabla c(t) \nabla\phi_1 dx + \frac{\alpha}{b} \int_{\partial\Omega}  \partial_t c_{\tau}(t)\phi_{1\tau} d\sigma + \frac{\alpha}{b}  \int_{\partial\Omega} \nabla_{\tau} c(t) \nabla_{\tau}\phi_1 d\sigma   \\
\   &\hspace{1cm}   +   \int_{\Omega} u\nabla c(t) \phi_1 dx= \int_{\Omega} -n(t) f(c(t))\phi_1 dx, \\
 \ &\int_{\Omega} \partial_t  n(t)\phi_2 dx +  \int_{\Omega} \left( \beta\nabla n(t) - g( n(t),\ c(t)) \nabla c(t)   \right)  \nabla\phi_2 dx \\
 \ &\hspace{3cm}     + \int_{\Omega} u(t)\nabla n(t) \phi_2 dx=0,\\
 \ &\int_{\Omega}\partial_t   u(t)\phi_3 dx +   \int_{\Omega} \xi  \nabla u(t) \nabla\phi_3 dx +  \int_{\Omega} (u(t)\cdot\nabla) u(t) \phi_3 dx=\\
 \ &\hspace{3cm}     \int_{\Omega} n(t) \nabla \sigma \phi_3 dx.
 \end{array} \right.
\end{equation}

for all $(\phi_1,\  \phi_2,\  \phi_3)\in (H^1(\Omega))^2\times V$ .

\end{theorem}

The paper is organized as follows.  In the next section, we introduce some preliminary lemmas and time-discretization scheme (Rothe's method).
We also outline the approaches to prove the main result.

Section 3 is devoted to proving Theorem 2.3.  We use regularization technique and fixed-point theorem to prove the existence of solution for
an auxiliary problem, then use Galerkin method to show existence of solutions to the discretized scheme.

In section 4, we derive several a priori estimates which will allow us to pass limits in the discretization scheme and thereby verify Theorem 1.1.

\section{Preliminaries }

In this paper, we will use the following Gronwall's lemma in the discrete form.

 \begin{lemma}
Let $0<k<1,\  (a_i)_{i\ge 1},\   $and $(A_i)_{i\ge 1}$ be sequence of real, non-negative numbers.  Assuming that  $(A_i)_{i\ge 1}$ 
is non-decreasing and that
\[  a_i \le A_i + k\sum_{j=0}^i  a_j, \quad \mbox{for}\  i=0,1,2,\cdots \]
then
\[  a_i \le \frac{1}{1-k} A_i exp \left(  (i-1)\frac{k}{1-k}   \right) ,  \quad \mbox{for}\  i=0,1,2,\cdots \]
\end{lemma}

 We also need the following relation:
 \begin{equation}\label{e2.1}
  2\int_{\Omega} a(a-b)dx = \Vert a\Vert _{L^2(\Omega)}^2 -  \Vert b\Vert _{L^2(\Omega)}^2 +   \Vert a-b \Vert _{L^2(\Omega)}^2.  
  \end{equation}
  
  \smallskip
  
 We will frequently use Young's inequality:
 \[ ab\le \delta a^2 + \frac{1}{4\delta} b^2.\]

 The compactness result presented in the next lemma will allow us to pass the limit in the Rothe approximation.
 \begin{lemma}
 Let $X,\ Y$ be two Banach spaces, such that $Y \subset X\ $, the injection being compact.
 
 Assume that $G$ is a family of functions in $L^{2}(0,\ T; \  Y)\bigcap L^p(0,\ T; \ X)$ for some $T>0$ and $p>1$, such that
 \begin{eqnarray}
 \    &   G \   \mbox{is bounded in}\   L^{2}(0,\ T; \  Y)\    \mbox{and}\   L^p(0,\ T; \ X);\nonumber\\
 \    &  \sup_{g\in G} \int_0^{T-a}  \Vert  g(a+s)- g(s)          \Vert _{X}^p ds \rightarrow 0  \quad  \mbox{as}\   a\rightarrow 0, \  a>0 \label{e2.2}
 \end{eqnarray}
 
 Then the family $G$ is relatively compact in $L^p(0,\ T; \ X)$.
 \end{lemma}

The proof the main result stated in section 1 is based on the Rothe's method of time discretization.  We divide the time interval 
$[0,\ T]$ into $N$ subintervals $[t_{m-1},\ t_m ]$, $t_m=km$, $k=T/N$, $m=1,\ 2,\ \cdots ,\ N.$

The time discretized variational formulation reads as

\begin{equation} \label {e2.3}
\left\{
\begin{array}{lcl}
\delta_t c^m -\alpha \Delta c^m  + u^m\cdot  \nabla c^m  + n^m f(c^m) =0   &\mbox{a.e.  in}  &\Omega\times(0,T), \\
\delta_t n^m  - \nabla\cdot  \left(\beta  \nabla n^m -g(n^m,\ c^m)  \nabla c^m \right)  + u^m \cdot \nabla n^m =0 &\mbox{a.e.  in}  &\Omega\times(0,T),\\
\delta_t u^m  - \nabla\cdot  (  \xi  \nabla u^m)  +(u^m\cdot \nabla) u^m +\nabla p^m = n^m \nabla   \sigma   &\mbox{a.e.  in}  &\Omega\times(0,T),\\
\nabla \cdot u^m =0 &\mbox{a.e.  in}  &\Omega\times(0,T),\\
\delta_t c^m =  \Delta_{\tau} c^m -b \partial_{\eta}c^m &\mbox{a.e.  on}  &\partial\Omega,\\
\beta  \partial_{\eta}n^m =  g(n^m,\ c^m)  \partial_{\eta}c^m &\mbox{a.e.  on}  &\partial\Omega,\\
u^m=0&\mbox{a.e.  on}  &\partial\Omega,\\
c^0=c_0(x),\   n^0=n_0(x),\   u^0=u_0(x).
\end{array} \right.
\end{equation}

Here, we use the notation:

$\delta_t c^m=\frac{c^m - c^{m-1}}{k},\  \delta_t n^m=\frac{n^m - n^{m-1}}{k}, \delta_t u^m=\frac{u^m - u^{m-1}}{k}$, where
$c^m, \ n^m, \ u^m, \  m=1,2,\cdots,N$ are the approximations of $c(x,t_m),\  n(x,t_m),\  $ and $u(x,t_m)$ respectively.

The existence result for the discrete scheme is stated in the following theorem, which will be proved in section 3.
\begin{theorem} \label{th2.3}
Suppose $(H_1)-(H_3)$ holds, $(c^0,\  n^0,\ u^0)\in (L^2(\Omega))^2\times H.$  Then there exists 
$(c^m,\  n^m,\ u^m)\in (L^2(\Omega))^2\times H $ solving the discrete problem \eqref{e2.3} for time step $k$ small enough.
\end{theorem}

With this result, we introduce the Rothe functions:
\begin{equation}
\left\{
\begin{array}{ll}
\  &\tilde{c}_k(t) = c^{m-1} + (t-t_{m-1})\delta_t c^m,\\  
\  &\tilde{c}_{k\tau}(t) = c_{\tau}^{m-1} + (t-t_{m-1})\delta_t c_{\tau}^m,\\  
\  &\tilde{n}_k(t) = n^{m-1} + (t-t_{m-1})\delta_t n^m,\\  
\  &\tilde{u}_k(t) = u^{m-1} + (t-t_{m-1})\delta_t u^m,\quad \mbox{for}\   t_{m-1}\le t\le t_m,\  1\le m\le N
\end{array} \right.
\end{equation}

and step functions:
\begin{equation}
\left\{
\begin{array}{ll}
\    & (c_k(t),\  c_{k\tau}(t),\ n_k(t),\  u_k(t))   =  (c^{m},\ c_{\tau}^{m},\  n^{m},\ u^{m}) ,\\  
\    & (c_k(0),\ c_{k\tau}(0),\ n_k(0),\  u_k(0))   =  (c_0,\ c_{0\tau},\ n_0,\ u_0).
\end{array}\right.
\end{equation}

for $t_{m-1}\le t\le t_m,\quad 1\le m \le N.$

We will prove the Rothe functions  $(\tilde{c}_k,\  \tilde{c}_{k\tau},\ \tilde{n}_k,\ \tilde{u}_k )$  and the step functions  $(c_k,\ c_{k\tau},\ n_k,\ u_k)$   will converge to the limit functions $(c,\  c_{\tau},\ n,\ u)$ as $k\rightarrow 0$,
and the limit functions will be the solutions to problem \eqref{e1.6}.  Thus, the main result will be proved. 

\section {Proof of Theorem 2.3}

The proof of Theorem 2.3 is based on semi-Galerkin method.  Let ${w_ j}$ be an orthogonal basis of $V$ that is orthonormal in $H$.
We denote by $V_j$ the finite vector space spanned by $ \{w_i\}_{1\le i \le j}.$

For a fixed $m$, let $u_j^m \in V_j$ be the Galerkin approximation of $u^m$, and assume $(c^{m-1},\ n^{m-1},\ u^{m-1})$ are given,
we consider the following problem:

Find $(c_j^{m},\ n_j^{m},\ u_j^{m})  \in  \left( H^1(\Omega)  \right)^2\times V_j$,  satisfying the following elliptic system:

\begin{equation}\label{e3.1}
\left\{
 \begin{array}{lll}
\  & -k\alpha\Delta c_j^m + k (u_j^m \cdot \nabla ) c_j^m + c_j^m = -k n_j^m f( c_j^m)+h 
 &\mbox{a.e.  in}\  \Omega\times(0,T),\\
\  &-k\nabla\cdot \left(  \beta\nabla n_j^m-g (n_j^m,\  c_j^m) \nabla c_j^m   \right) + k (u_j^m \cdot \nabla ) n_j^m + n_j^m = l  &\mbox{a.e.  in}\  \Omega\times(0,T),\\
\  &-k\nabla\cdot (\xi\nabla u_j^m) + k (u_j^m \cdot \nabla ) u_j^m + k \nabla P_j^m + u_j^m=k n_j^m\nabla \sigma + q &\mbox{a.e.  in}\  \Omega\times(0,T),\\
\ &\nabla\cdot u_j^m = 0 &\mbox{a.e.  in}\  \Omega\times(0,T),\\
 \  & k\partial_{\eta} c_j^m =  \frac{k}{b}   \Delta_{\tau} c_{j_\tau}^m -\frac{1}{b} c_{j_\tau}^m+  \frac{1}{b} h_{\tau}  &\mbox{a.e.  on}\  \partial\Omega, \\
 \  & \beta  \partial_{\eta}n_j^m =  g(n_j^m,\ c_j^m)  \partial_{\eta}c_j^m&\mbox{a.e.  on}\  \partial\Omega,\\
 \ &u_j^m=0 &\mbox{a.e.  on}\  \partial\Omega.
 \end{array}
\right.
\end{equation}

where $(h,\ h_{\tau},\ l,\ q)=(c^{m-1},\ c_{\tau}^{m-1},\ n^{m-1},\ u^{m-1}).$

 The existence result of the above problem is stated in the following theorem:
\begin{theorem}\label{th3.1}
Suppose $(H_1)-(H_3)$ holds, and $(h,\  l,\ q,\ h_{\tau})\in (L^2(\Omega))^3\times L^2(\partial\Omega),$  then there exists a solution
$(c_j^m,\  n_j^m,\ u_j^m)\in (H^1(\Omega))^2\times V $ of problem \eqref{e3.1}.
\end{theorem}

 We will then move on to derive estimates on the solution $(c_j^m,\ c_j^m,\ c_j^m)$ of problem \eqref{e3.1}.  Those estimates allow us to pass the limits
 by letting $j\rightarrow\infty$, the sequences will converge to limit functions, which will be the solutions of the discrete problem \eqref{e2.3}, 
 and Theorem 2.3 will be proved. 

\smallskip

The proof of Theorem 3.1 will use the following lemma:
\begin{lemma}
Suppose $(H_1)-(H_3)$ holds, and $(h,\  l,\ h_{\tau})\in (L^2(\Omega))^2\times L^2(\partial\Omega),$ then for fixed $\hat{u} \in V_j$, there exist a 
solution $(c,\ n)\in (H^1(\Omega))^2$ of the following problem:
\begin{equation}\label{e3.2}
\left\{
 \begin{array}{ll}
\  & -k\alpha\Delta c + k (\hat{u}\cdot \nabla  c )+ c = -k n f( c)+h,\\
\  &-k\nabla\cdot \left(  \beta\nabla n-g (n,\  c) \nabla c   \right) + k (\hat{u}\cdot \nabla  n) + n =  l,\\
 \  & k\partial_{\eta} c =  \frac{k}{b}   \Delta_{\tau} c_{\tau} -\frac{1}{b} c_{\tau}+  \frac{1}{b} h_{\tau},  \\
 \  & \beta  \partial_{\eta}n =  g(n,\ c)  \partial_{\eta}c.
 \end{array}
\right.
\end{equation}

\end{lemma}

\bigskip

3.1  Proof of Lemma 3.2

\smallskip

To prove Lemma 3.2, we first consider a regularized problem, we smooth $\hat{u}$ by replacing  $\hat{u}$ with $J_{\epsilon}\hat{u}$,
where $J_{\epsilon}\hat{u} = \left(   (\psi _{\epsilon} u )\ast w_{\epsilon}   \right)_{div}.$

Here 

\[
\psi _{\epsilon}(x):= \left\{ \begin{array}{ll}
        0 & \mbox{if}\   \   \mbox{dis}\  (x,\ \partial \Omega) \le 2\epsilon,\\
         1  & \mbox{elsewhere} 
                \end{array}\right.
 \]

$(\psi _{\epsilon} u )\ast w_{\epsilon}$ denote the standard regularization of $\psi _{\epsilon} u $ with kernel $w_{\epsilon}$ having support
in a ball of radius $\epsilon$.  The symbol $(\cdot)_{div}$ comes from the Helmholtz decomposition, see \cite{bul}.

$J_{\epsilon}\hat{u}$ keeps Dirichlet boundary condition and divergence free property, therefore, we have the identity
$\int_{\Omega}  (J_{\epsilon} \hat{u}\cdot \nabla  c) c dx=0$, see \cite{tem}.  We will use this identity frequently throughout this paper.

We have the following existence result for the regularized problem.
\begin{lemma}
Suppose $(H_1)-(H_3)$ holds, and $(h,\  l,\ h_{\tau})\in (L^2(\Omega))^2\times L^2(\partial\Omega),$ then for fixed $\epsilon\in (0,\ 1)\ $and $\hat{u} \in V_j$, there exist a 
solution $(c_{\epsilon},\ n_{\epsilon})\in (H^1(\Omega))^2$ of the following problem:
\begin{equation}\label{e3.3}
\begin{array}{ll}
\  & -k\alpha\Delta c_{\epsilon} + k J_{\epsilon}\hat{u}\cdot \nabla  c_{\epsilon} + c_{\epsilon} = -k n_{\epsilon} f( c_{\epsilon})+h,  \\
\  &-k\nabla\cdot \left(  \beta\nabla n_{\epsilon}-g (n_{\epsilon},\  c_{\epsilon}) \nabla c_{\epsilon}   \right) + k J_{\epsilon}\hat{u}\cdot \nabla  n_{\epsilon} + n_{\epsilon} =   l,  \\
 \  & k\partial_{\eta} c_{\epsilon} =  \frac{k}{b}   \Delta_{\tau} {c_{\epsilon}}_{\tau} -\frac{1}{b} {c_{\epsilon}}_{\tau}+  \frac{1}{b} h_{\tau},  \\
 \  & \beta  \partial_{\eta}n_{\epsilon} =  g(n_{\epsilon},\ c_{\epsilon})  \partial_{\eta}c_{\epsilon}.
 \end{array}
\end{equation}
\end{lemma}

{\bf Proof of Lemma 3.3}:  In the proof, we omit subscript $\epsilon$, write $(c_{\epsilon},\  n_{\epsilon})$ as $(c,\ n)$.  We use Schaefer's 
fixed-point theorem.

We define an operator $\Phi :\quad X \rightarrow X,\ $ where 
\[ X= H^1(\Omega) \times H^1(\Omega)  \]

Fix $(\hat{c}, \   \hat{n})\  \in X,\ $  we set $\Phi (\hat{c}, \   \hat{n}) = (c,\ n)$,
 where $(c,\ n)$ is the solution to the following problem:
 \begin{eqnarray}
&\ & -k\alpha\Delta c + k J_{\epsilon} \hat{u}\cdot \nabla  c + c = -k \hat{n} f( \hat{c})+h,\label{e3.4}
\\
&\ & -k\nabla\cdot \left(  \beta\nabla n - g (\hat{n},\  c) \nabla c   \right)+ k J_{\epsilon} \hat{u}\cdot \nabla  n + n =    l,\label{e3.5}\\
&\ & k\partial_{\eta} c =  \frac{k}{b}   \Delta_{\tau} c_{\tau} -\frac{1}{b} c_{\tau}+  \frac{1}{b} h_{\tau}, \label{e3.6} \\
&\ & \beta  \partial_{\eta}n =  g(\hat{n},\ c)  \partial_{\eta}c. \label{e3.7}
 \end{eqnarray}

 We want to show $\Phi$ is a continuous and compact mapping of $X$ into itself, such that the set 
 $\{ x\in X:\   x=\lambda \Phi (x)\   \mbox{for some}\   0\le\lambda\le 1\} $ is bounded.  Therefore, by
 Schaefer's fixed-point theorem, $\Phi$ has a fixed point.

 Multiply both sides of the equation \eqref{e3.4} by $c$, use $\int_{\Omega}  (J_{\epsilon} \hat{u}\cdot \nabla  c) c dx=0$,     
 we arrive at
\[
  \begin{array}{ll}
 \   &  -\alpha\int_{\partial \Omega} c ( h_2 -  \frac{1}{b} c_{\tau} +  \frac{k}{b}   \Delta_{\tau} c_{\tau} ) d\sigma 
  +\alpha k \Vert \nabla c \Vert _{L^2}^2  +  \Vert  c \Vert _{L^2}^2 \\
 \   &\hspace{1cm}   =   \int_{\Omega} h_1 c.
 \end{array}
\]

\smallskip

Here $ h_1=  -k\hat{n}f(\hat{c}) + h    ,\     h_2 = \frac{1}{b} h_{\tau}.$

\smallskip

Use H$\ddot{o}$lder's inequality and Young's inequality, we get

\begin{equation}\label{e3.8}
  \begin{array}{ll}
 \   &  \frac{\alpha k}{b} \Vert \nabla c_{\tau} \Vert _{L^2}^2  +  {\alpha k} \Vert \nabla c \Vert _{L^2}^2 +  (1-\delta) \Vert  c \Vert _{L^2}^2
 +\alpha (\frac{1}{b}-\delta)\Vert  c_{\tau} \Vert _{L^2}^2 \\
 \   &\hspace{0.1cm}   \le  \frac{1}{4\delta} \Vert h_1\Vert_{L^2}^2  +   \frac{\alpha}{4\delta} \Vert h_2\Vert_{L^2}^2\\
  \   &\hspace{0.1cm}   \le  \frac{1}{2\delta} \Vert h\Vert_{L^2}^2  +  \frac{\alpha}{4\delta b^2 } \Vert h_{\tau}\Vert_{L^2}^2
                                        +  \frac{1}{2\delta} k^2 f_1^2 \Vert \hat{n}\Vert_{L^2}^2.
 \end{array}
\end{equation}

Choose $\delta < \min (1,\frac{1}{b})$ in the above inequality, we have 
\begin{equation}\label{e3.8-1}
\Vert  c \Vert _{L^2}^2 +  \Vert \nabla c \Vert_{L^2}^2 \le C(\Vert h\Vert_{L^2}^2 + \Vert h_{\tau}\Vert_{L^2}^2 + \Vert \hat{n}\Vert_{L^2}^2).
\end{equation}

\smallskip

Multiply \eqref{e3.5} by $n$, integrate by parts, and use boundary conditions, we get,
\[
\begin{array}{ll}
 \   & k\beta\Vert \nabla n \Vert _{L^2}^2 + \Vert n \Vert _{L^2}^2  =   \int_{\Omega} ln + k\int_{\Omega}  g (\hat{n},\  c)  \nabla c \nabla n     \\
\   &\hspace{1cm} \le  \delta \Vert n \Vert _{L^2}^2  + \frac{1}{4\delta} \Vert l \Vert _{L^2}^2 + \frac{kg_1}{2} (\Vert \nabla n \Vert _{L^2}^2
+\Vert \nabla c \Vert _{L^2}^2),
 \end{array}
\]
which gives,
\[
(k\beta-\frac{kg_1}{2} )\Vert \nabla n \Vert _{L^2}^2 + (1-\delta) \Vert n \Vert _{L^2}^2 \le  \frac{1}{4\delta} \Vert l \Vert _{L^2}^2
+\frac{kg_1}{2}\Vert \nabla c \Vert _{L^2}^2.
\]
Choose $\delta < 1$ in the above inequality, and require $\beta > \frac{g_1}{2}$, we have
\begin{equation}\label{e3.9-1}
\Vert  n \Vert _{L^2}^2 +  \Vert \nabla n \Vert_{L^2}^2 \le C(\Vert l\Vert_{L^2}^2 +\Vert \nabla c \Vert_{L^2}^2 ).
\end{equation}
This together with \eqref{e3.8-1} shows
\begin{equation}\label{e3.9-2}
\Vert  n \Vert _{L^2}^2 +  \Vert \nabla n \Vert_{L^2}^2 \le C(\Vert l\Vert_{L^2}^2 +\Vert h\Vert_{L^2}^2 + \Vert h_{\tau}\Vert_{L^2}^2 + \Vert \hat{n}\Vert_{L^2}^2 ).
\end{equation}

\eqref{e3.8-1} and \eqref{e3.9-2} show that $\Phi$ maps $X$ to $X$.

Next, we want to show the set 
\[ S= \{ (c,\ n)\in H^1 \times H^1 :\   (c,\ n)=\lambda \Phi (c,\ n)\   \mbox{for some}\   0\le\lambda\le 1\}  \]

is bounded.

\smallskip

We consider the equation $ (\frac{c}{\lambda},\   \frac{n}{\lambda} )= \Phi (c,\ n)$:
 \begin{eqnarray}
&\ & -k\alpha\Delta c + k J_{\epsilon} \hat{u}\cdot \nabla  c + c = \lambda ( -k n f( c)+h ),\label{e3.15}\\
&\ &-k\beta\Delta n + k J_{\epsilon} \hat{u}\cdot \nabla  n + n =  \lambda\left( -k   \nabla\cdot  (g (n,\  \frac{c}{\lambda}) 
\frac{\nabla c}{ \lambda})  +  l \right),\label{e3.16}\\
&\ & k\partial_{\eta} c =  \frac{k}{b}   \Delta_{\tau} c_{\tau} -\frac{1}{b} c_{\tau}+  \frac{\lambda}{b} h_{\tau}, \label{e3.17} \\
 &\ & \beta  \partial_{\eta}n =  g(n,\   \frac{c}{\lambda})  \partial_{\eta}c.\label{e3.18}
 \end{eqnarray}

Multiply \eqref{e3.15} by $c$, integrate with respect to $x$, we get
\[
 \begin{array}{ll}
 \   & \frac{\alpha k}{b}\Vert  \nabla_{\tau} c_{\tau} \Vert _{L^2}^2 + \frac{\alpha }{b}\Vert  c_{\tau} \Vert _{L^2}^2
- \frac{\alpha \lambda}{b}\int_{\partial\Omega}ch_{\tau}d\sigma + k\alpha\Vert  \nabla c \Vert _{L^2}^2
+\Vert  c \Vert _{L^2}^2\\
\   &\hspace{1cm}   = -\lambda\int_{\Omega}knf(c)c +  \lambda\int_{\Omega} hc.
 \end{array}
 \]
 
 Use Young's inequality, we arrive at
\[
 \begin{array}{ll}
 \   & \frac{\alpha k}{b}\Vert  \nabla_{\tau} c_{\tau} \Vert _{L^2}^2 + ( \frac{\alpha }{b} - \frac{\alpha \lambda}{b}\delta)
 \Vert  c_{\tau} \Vert _{L^2}^2 + k\alpha\Vert  \nabla c \Vert _{L^2}^2 + (1-\lambda k f_1 \delta - \lambda\delta)
 \Vert  c \Vert _{L^2}^2\\
 \   &\hspace{1cm} \le  \frac{\alpha \lambda}{4b\delta}\Vert  h_{\tau} \Vert _{L^2}^2 + \frac{\lambda k f_1}{4 \delta}
 \Vert  n \Vert _{L^2}^2 + \frac{\lambda}{4 \delta} \Vert  h \Vert _{L^2}^2.
 \end{array}
 \]
 Choose $\delta < \min(\frac{1}{\lambda}, \frac{1}{\lambda +\lambda k f_1 })$ in the above inequality, we have
\begin{equation}\label{e3.19}
 \begin{array}{ll}
 \   &   \alpha k \Vert  \nabla c \Vert _{L^2}^2  
\le  \frac{1}{4\delta}\left( {\lambda}  k f_1 \Vert n\Vert_{L^2}^2    
  + {\lambda}  \Vert h\Vert_{L^2}^2\right) + \frac{ {\lambda}\alpha}{4\delta b } \Vert h_{\tau}\Vert_{L^2}^2.
 \end{array}
 \end{equation}
 
 Multiply \eqref{e3.16} by $n$, integrate with respect to $x$, use $ \int_{\Omega} (J_{\epsilon} \hat{u})\cdot (\nabla  n) n = 0$,
then integrate by parts, use boundary condition \eqref{e3.18}, H$\ddot{o}$lder and Young's inequalities, we arrive at:
\[
k\beta\Vert  \nabla n \Vert _{L^2}^2 + \Vert  n \Vert _{L^2}^2 \le  \lambda (\hat{\delta} \Vert  n \Vert _{L^2}^2 + 
\frac{1}{ 4{ \hat{\delta} }  } \Vert  l \Vert _{L^2}^2)  +  kg_1( \frac{1}{ 4{ \hat{\delta} }  } \Vert  \nabla c \Vert _{L^2}^2
+ \hat{\delta} \Vert  \nabla n \Vert _{L^2}^2),
\]
which gives,
\begin{equation}\label{e3.20}
\begin{array}{ll}
\   & (1-\lambda \hat{\delta}) \Vert  n \Vert _{L^2}^2 + (k\beta - k g_1 \hat{\delta}) \Vert  \nabla n \Vert _{L^2}^2 \\
\   & \le \frac{\lambda}{4\hat{\delta}} \Vert  l \Vert _{L^2}^2 +  \frac{k g_1}{4 \hat{\delta}} \Vert  \nabla c \Vert _{L^2}^2.
 \end{array}
\end{equation}

Choose $\hat{\delta} <  \min (\frac{1}{\lambda},\  \frac{\beta}{g_1} )$ in \eqref{e3.20}, with fixed $\delta,\ \hat{\delta}$ in \eqref{e3.19} and \eqref{e3.20}, choose $k$ small, such that
\[  \frac{ {\lambda}  k f_1}{4\delta} < 1-\lambda \hat{\delta}. \]

Require $\alpha$ to be a constant s.t. $\alpha > \frac{g_1}{4\hat{\delta}}.$  By adding \eqref{e3.19} and \eqref{e3.20}, we can absorb the terms with
$\Vert  \nabla c \Vert _{L^2}^2$ and $\Vert  n \Vert _{L^2}^2$ in the right hand side of the equation to the left hand side, 
and obtain:
\begin{equation}\label{e3.21}
\begin{array}{ll}
 \   &  (\alpha k-  \frac{k g_1}{4 \hat{\delta}})  \Vert  c \Vert _{H^1}^2 
 +  (k\beta - k g_1 \hat{\delta}) \Vert  \nabla n \Vert _{L^2}^2 + (  1-  \lambda \hat{\delta} -  \frac{ {\lambda}  k f_1 }
 {4\delta} )  \Vert n\Vert_{L^2}^2 \\
\  & \le   \frac { \lambda}{ 4\delta }  \Vert h\Vert_{L^2}^2 +  \frac{ {\lambda}\alpha}{4\delta b} \Vert \ h_{\tau} \Vert_{L^2}^2
+  \frac { \lambda}{ 4\hat{\delta} } \Vert l\Vert_{L^2}^2.
\end{array}
\end{equation}
Therefore, $\Vert  n \Vert _{H^1}, \   \Vert  c \Vert _{H^1}$ is bounded, hence the set $S$ is bounded.

Next, we want to check the mapping $\Phi$ is compact. 

We proceed to show that $\Phi$ maps bounded set
 $(\hat{c}, \   \hat{n}) \in H^1 \times H^1 $ to bounded set  $(c,\  n) = \Phi(\hat{c}, \   \hat{n}) $ in $H^2 \times H^2 $:
 
 For $(\hat{c}, \   \hat{n})$  bounded in $H^1 \times H^1 $, from equation \eqref{e3.4},  \eqref{e3.6}, and elliptic regularity theorem, we know that
$\Vert c\Vert_{H^2}$ is bounded.  Boundedness of $\Vert n\Vert_{H^2}$ comes from  equation  \eqref{e3.5},  \eqref{e3.7}, elliptic regularity theorem, 
and boundedness of $\Vert c\Vert_{H^2}$.  So we have $(c,\  n) = \Phi(\hat{c}, \   \hat{n}) $ is bounded in 
$H^2 \times H^2 $.

Since the embedding of $H^2 \times H^2 $ into $H^1 \times H^1 $ is compact, the bounded set  $(c,\  n)$
in $H^2 \times H^2 $ is relatively compact in $H^1 \times H^1 $.  Therefore, the mapping $\Phi$ is compact.

The last step is to show that $\Phi$ is continuous.

Let $(\hat{c}_n,\  \hat{n}_n)$ converge to $(\hat{c},\  \hat{n})$ in $H^1 \times H^1$ strongly.  Let $(c_n,\  n_n)
=\Phi (\hat{c}_n,\  \hat{n}_n)$, so we have
\begin{equation}\label{e3.22}
 \begin{array}{ll}
\  & -k\alpha\Delta c_n + k J_{\epsilon} \hat{u}\cdot \nabla  c_n + c_n =   -k\hat{n}_n f( \hat{c}_n)+h, \\
\  &-k\beta\Delta n_n + k J_{\epsilon} \hat{u}\cdot \nabla  n_n + n_n =   -k   \nabla\cdot  (g (\hat{n}_n,\  c_n) 
\nabla c_n)  +  l ,\\
 \  & k\partial_{\eta} c_n =  \frac{k}{b}   \Delta_{\tau} {c_n}_{\tau} -\frac{1}{b} {c_n}_{\tau}+  \frac{1}{b} h_{\tau},  \\
 \  & \beta  \partial_{\eta}n_n =  g  (\hat{n}_n,\  c_n)\partial_{\eta}c_n.
 \end{array}
\end{equation}

Since the sequence $(\hat{c}_n,\  \hat{n}_n)$ is bounded in  $H^1 \times H^1$, the same argument to show 
the compactness of $\Phi$ can be applied here to derive boundedness of $(c_n,\  n_n)$ in 
$H^2 \times H^2$.  Since $H^2 \times H^2$ is compactly embedded into   $H^1 \times H^1$, there exist
$(c,\  n) \in  H^2 \times H^2$ and subsequence of $(c_n,\  n_n)$, still denoted as  $(c_n,\  n_n)$, s.t.
\begin{eqnarray}
&\  & (c_n,\  n_n) \rightharpoonup (c,\  n) \    \mbox{weakly in} \    H^2 \times H^2,  \label{e3.23}\\
&\  & (c_n,\  n_n) \rightarrow (c,\  n) \    \mbox{strongly in} \    H^1 \times H^1.  \label{e3.24}
\end{eqnarray}

With \eqref{e3.24} and the fact that $(\hat{c}_n,\  \hat{n}_n)$ converges to $(\hat{c},\  \hat{n})$ in $H^1 \times H^1$ strongly,
we have
\begin{equation} \label{e3.25}
\begin{array}{ll}
\  &(\hat{c}_n,\  \hat{n}_n,\  c_n,\  n_n) \rightarrow (\hat{c},\  \hat{n},\ c,\  n) \    \mbox{a.e. in} \   \Omega.
\end{array}
\end{equation}

As for the trace of $(c_n,\  n_n)$,  we have $ {n_n}_{\tau}$ is bounded in $H^{ \frac{3}{2} } (\partial \Omega)\ $ (since 
$n_n$ is bounded in $H^2$), therefore, is also bounded in $H^1(\partial \Omega)$.  Since $H^1(\partial \Omega)\ $ is compactly embedded into $L^2(\partial \Omega)$, there exists subsequence of $ {n_n}_{\tau}$, still denoted as $ {n_n}_{\tau}$, s.t.
$ {n_n}_{\tau} \rightarrow {n}_{\tau}$ strongly in $ L^2(\partial \Omega)$, so we have
 \begin{equation} \label{e3.26}
 {n_n}_{\tau} \rightarrow {n}_{\tau}\   \mbox{almost everywhere}.  
 \end{equation}

 Similarly, we have  
 \begin{equation} \label{e3.27}
  {c_n}_{\tau} \rightarrow {c}_{\tau} \   \mbox{almost everywhere}.
 \end{equation}

With assumption $(H_1) - (H_3)$,  \eqref{e3.25} -  \eqref{e3.27},  and dominated convergence theorem, we  conclude that
\begin{eqnarray}
&\  & f(\hat{c}_n) \rightarrow  f(\hat{c}) \    \mbox{a.e.  and strongly in} \    L^2, \label{e3.28}\\
&\  & g(\hat{n}_n,\  c_n) \rightarrow g(\hat{n},\  c) \    \mbox{a.e.  and strongly in} \    L^2. \label{e3.29}
\end{eqnarray}

From boundedness of $ (c_n,\  n_n)$ in $H^1 \times H^1$, we have
\begin{equation}\label{e3.30}
(\nabla c_n,\  \nabla n_n) \rightharpoonup (\nabla c,\  \nabla n) \    \mbox{weakly in} \   L^2\times L^2 .
\end{equation}

From \eqref{e3.28} -  \eqref{e3.30}, we get
\[
\begin{array}{ll}
\  &\hat{n}_n f(\hat{c}_n) \rightarrow \hat{n} f(\hat{c}) \    \mbox{in distribution,}  \\
\  & g(\hat{n}_n,\  c_n) \nabla c_n \rightarrow g(\hat{n},\  c) \nabla c\    \mbox{in distribution.} 
\end{array}
\]

Similarly, we can pass limits in other terms in  \eqref{e3.22} by letting $n\rightarrow \infty$, we then obtain  \eqref{e3.4} -  \eqref{e3.7}, therefore,
we have shown that $(c,\ n) = \Phi (\hat{c},\  \hat{n}).$  The continuity of $\Phi$ is proved.   This finishes the proof
of Lemma 3.3.

{\bf Proof of Lemma 3.2:}  We multiply $\eqref{e3.3}_1$ by $c_{\epsilon}$ and $\eqref{e3.3}_2$ by $n_{\epsilon}$, integrate over $\Omega$, then proceed in an analogous way as in \eqref{e3.19} - \eqref{e3.21}, we obtain the following estimates:
\begin{equation}\label{e3.31}
\begin{array}{ll}
\ &\Vert c_{\epsilon\tau}   \Vert_{H^1}^2 + \Vert c_{\epsilon}   \Vert_{H^1}^2 + \Vert n_{\epsilon}   \Vert_{H^1}^2
+\Vert n_{\epsilon}   \Vert_{L^2}^2 \le C \left(  \Vert h  \Vert_{L^2}^2  + \Vert h_{\tau}  \Vert_{L^2}^2  
+\Vert l  \Vert_{L^2}^2 \right).
\end{array}
\end{equation}

where the constant $C$ is independent of $\epsilon$.

We obtain from \eqref{e3.31} that $(c_{\epsilon\tau},\  c_{\epsilon},\  n_{\epsilon})$ is bounded in 
$H^1(\partial \Omega) \times H^1( \Omega)   \times  H^1(\Omega),$  uniformly w.r.t.  $ \epsilon$.  Since
$H^1(\partial \Omega)$ and $H^1(\Omega)$ are compactly embedded in $L^2(\partial \Omega)$ and $L^2(\Omega)$
respectively, there exist $(c_{\tau},\  c,\  n) \in H^1(\partial \Omega) \times H^1( \Omega)   \times  H^1(\Omega) $
and subsequence of $(c_{\epsilon\tau},\  c_{\epsilon},\  n_{\epsilon})$, still denoted as $(c_{\epsilon\tau},\  c_{\epsilon},\  n_{\epsilon})$, s.t. as $\epsilon \rightarrow 0^{+}$, we have
\begin{equation}\label{e3.32}
\begin{array}{ll}
\  & (c_{\epsilon\tau},\  c_{\epsilon},\  n_{\epsilon}) \rightharpoonup (c_{\tau},\  c,\  n) \    \mbox{weakly in} \ 
H^1(\partial \Omega)\times H^1(\Omega) \times H^1(\Omega),\\
\  & (c_{\epsilon\tau},\  c_{\epsilon},\  n_{\epsilon}) \rightarrow (c_{\tau},\  c,\  n) \    \mbox{strongly in} \ 
L^2(\partial \Omega)\times L^2(\Omega) \times L^2(\Omega).
\end{array}
\end{equation}

From \eqref{e3.32} together with 
\begin{equation}\label{e3.33}
\begin{array}{ll}
\   & J_{\epsilon} \hat{u} \rightarrow  \hat{u} \     \mbox{strongly in} \   L^2(\Omega).
\end{array}
\end{equation}

 We get
\[
\begin{array}{ll}
\   & J_{\epsilon} \hat{u} \cdot \nabla c_{\epsilon}  \rightarrow  \hat{u}\cdot \nabla c\    \mbox{in distribution,} \\   
\   & J_{\epsilon} \hat{u} \cdot \nabla n_{\epsilon}  \rightarrow  \hat{u}\cdot \nabla n\    \mbox{in distribution.} \\   
\end{array}
\]

We then proceed in an analogous way as in the proof of continuity of $\Phi$ in Lemma 3.3 to pass the limit in
\eqref{e3.3} by letting $\epsilon \rightarrow 0^{+}$ to obtain \eqref{e3.2}.  The proof of Lemma 3.2 is completed.

\bigskip

3.2   Proof of Theorem 3.1

\smallskip

For the simplicity of notation, in the proof of Theorem 3.1, we omit the superscript ``$m$", and write equation \eqref{e3.1} as follows:
\begin{equation}\label{e3.34}
\left\{
 \begin{array}{lll}
\  & -k\alpha\Delta c_j + k (u_j \cdot \nabla ) c_j + c_j = -k n_j f( c_j)+h\  &\mbox{in}\  \Omega,\\
\  &-k\nabla\cdot \left(  \beta\nabla n_j-g (n_j,\  c_j) \nabla c_j   \right) + k (u_j \cdot \nabla ) n_j + n_j = l\  &\mbox{in}\  \Omega,\\
\  &-k\nabla\cdot (\xi\nabla u_j) + k (u_j \cdot \nabla ) u_j + k \nabla P_j + u_j=k n_j\nabla \sigma + q\  &\mbox{in}\  \Omega,\\
\ &\nabla\cdot u_j = 0\  &\mbox{in}\  \Omega,\\
 \  & k\partial_{\eta} c_j =  \frac{k}{b}   \Delta_{\tau} c_{j_\tau} -\frac{1}{b} c_{j_\tau}+  \frac{1}{b} h_{\tau}\  &\mbox{on}\  \partial\Omega,  \\
 \  & \beta  \partial_{\eta}n_j =  g(n_j,\ c_j)  \partial_{\eta}c_j\  &\mbox{on}\  \partial\Omega,\\
 \ &u_j=0 \  &\mbox{on}\  \partial\Omega.
\end{array}\right.
\end{equation}

Recall that $u_j\in V_j$ is the Galerkin approximation of $u^m$ for a fixed $m$.  We define an operator
$L:\  \    V_j \rightarrow V_j$ as follows:

For a fixed $\hat{u}_j\in V_j,\ $ we set  $L(\hat{u}_j)=u_j,\  $ where $u_j$ is the solution to the following problem:
 \begin{eqnarray}
&\  & -k\alpha\Delta c_j + k (\hat{u}_j \cdot \nabla ) c_j + c_j = -k n_j f( c_j)+h, \label{e3.35}\\
&\  &-k\nabla\cdot \left(  \beta\nabla n_j - g (n_j,\  c_j) \nabla c_j   \right) + k (\hat{u}_j\cdot \nabla ) n_j + n_j = l, \label{e3.36}\\
&\  &-k\nabla\cdot (\xi\nabla u_j) + k (u_j \cdot \nabla ) u_j + k \nabla P_j + u_j=k n_j\nabla \sigma + q, \label{e3.37}\\
&\ &\nabla\cdot u_j = 0.\label{e3.38}
\end{eqnarray}

with boundary conditions:
 \begin{eqnarray}
&\  & k\partial_{\eta} c_j =  \frac{k}{b}   \Delta_{\tau} c_{j_\tau} -\frac{1}{b} c_{j_\tau}+  \frac{1}{b} h_{\tau},  \label{e3.39} \\
 &\  & \beta  \partial_{\eta}n_j =  g(n_j,\ c_j)  \partial_{\eta}c_j, \label{e3.40} \\
 &\ &u_j=0.  \label{e3.41} 
\end{eqnarray}

Let $\hat{u}_j$ belongs to a bounded set in $V_j$, we fix $\hat{u}_j$ in \eqref{e3.35} -  \eqref{e3.36},  solve  \eqref{e3.35}, \eqref{e3.36}, \eqref{e3.39}, \eqref{e3.40} for 
$(c_j,\  n_j)$.  The existence of solution $(c_j,\  n_j)$ is proved in Lemma 3.2, and we have the estimate
of $n_j$:
\begin{equation}\label{e3.42}
\Vert n_j \Vert_{L^2}^2 \le C ( \Vert h \Vert_{L^2}^2 + \Vert h_{\tau} \Vert_{L^2}^2 +\Vert l \Vert_{L^2}^2).
\end{equation}

We then use this $n_j$ in equation \eqref{e3.37}, proceed to solve equation \eqref{e3.37},\eqref{e3.38} and \eqref{e3.41} for $u_j$, since $u_j\in V_j$, which is finite-dimensional, the existence of this $u_j$ can be proved in an analogous way as in \cite[p. 164]{tem}.  The operator $L$ then maps $\hat{u}_j$ to $u_j$, i.e. $u_j = L(\hat{u}_j).$
 
 Next, we show that $L$ has a fixed point.
 
 To this end, we multiply equation \eqref{e3.37} with $u_j$, integrate over $\Omega$, use 
 $\int_{\Omega} (u_j\cdot \nabla )u_j\cdot u_j = 0,\ $ we get,
\[
\begin{array}{ll}
\  & k\xi \Vert \nabla u_j \Vert_{L^2}^2 + \Vert u_j \Vert_{L^2}^2 = k \int_{\Omega} n_j \nabla \sigma u_j 
+ \int_{\Omega} q  u_j. 
 \end{array}
 \]

Use Young's inequality for the terms on the right hand side of the above equation, and use \eqref{e3.42},  we obtain,
\begin{equation}\label{e3.43}
\begin{array}{ll}
\  & k\xi \Vert \nabla u_j \Vert_{L^2}^2 + ( 1-kC\delta- \delta )\Vert u_j \Vert_{L^2}^2\\
\  &\hspace{1cm} \le  \frac{kC}{4\delta} \Vert n_j \Vert_{L^2}^2 + \frac{1}{4\delta} \Vert q \Vert_{L^2}^2\\
\  &\hspace{1cm} \le  \frac{kC}{4\delta}( \Vert h \Vert_{L^2}^2 + \Vert h_{\tau} \Vert_{L^2}^2 +\Vert l \Vert_{L^2}^2)
+ \frac{1}{4\delta} \Vert q \Vert_{L^2}^2.
 \end{array}
\end{equation}

By choosing $\delta$ small enough, from the above inequality, we see that $ u_j = L(\hat{u}_j) $
is bounded in $V_j$. Therefore, $L$ maps a bounded set in $V_j$ to a bounded set in $V_j$.  
Since $V_j$ is finite dimensional space, we can use Brouwer fixed-point theorem to conclude that $L$ has a fixed point, which is the solution
of \eqref{e3.34}.  This completes the proof of Theorem 3.1.

\bigskip

3.3   Proof of Theorem 2.3

\smallskip

With Theorem 3.1, we have shown the existence results for \eqref{e3.1} with a fixed ``$m$".

Now we multiply $\eqref{e3.34}_1$ by $c_j$,  $\eqref{e3.34}_2$ by $n_j$,   $\eqref{e3.34}_3$ by $u_j$, integrate over $\Omega$, proceed in an analogous way as
in  \eqref{e3.19} - \eqref{e3.21}, \eqref{e3.43}, we see that $(c_j,\ n_j,\ u_j)$ of \eqref{e3.34} is bounded in $\left(H^1(\Omega)\right)^2 \times V$, uniformly 
with respect to $j$.

As a result, as $j\rightarrow +\infty$, we have
\[
\begin{array}{ll}
\  & (c_j,\  n_j) \rightharpoonup (c,\  n) \    \mbox{weakly in} \    H^1 \times H^1,\\
\  & (c_j,\  n_j) \rightarrow (c,\  n) \    \mbox{strongly in} \    L^2 \times L^2,\\
\  & u_j \rightharpoonup u \    \mbox{weakly in} \    V,\\
\  & u_j \rightarrow u \    \mbox{strongly in} \    H.
\end{array}
\]

We can pass the limit (letting $j\rightarrow +\infty$) in \eqref{e3.34}.  This finish the proof of Theorem 2.3.

\section{Proof of Theorem 1.1}

We first derive aprori estimates for $(c^m,\ n^m, \  u^m)$.  To this end, we multiply $\eqref{e2.3}_1$ by $c^m$,
$\eqref{e2.3}_2$ by $n^m$, $\eqref{e2.3}_3$  by $u^m$, and integrate over $\Omega$, we arrive at,
\[
\begin{array}{ll}
 \   &  \frac{1}{k}\int_{\Omega}   c^m (c^m-c^{m-1}) - \alpha \int_{\Omega} c^m  \Delta c^m  +
 \int_{\Omega} c^m n^m f(c^m) = 0,\\
  \   &  \frac{1}{k}\int_{\Omega}   n^m (n^m-n^{m-1}) - \int_{\Omega} n^m \nabla \cdot\left(\beta \nabla n^m -g(n^m,\ c^m)  \nabla c^m \right)=0,\\
  \   &  \frac{1}{k}\int_{\Omega}   u^m (u^m-u^{m-1}) - \int_{\Omega} u^m \nabla \cdot ( \xi  \nabla u^m)=
  \int_{\Omega} u^m n^m \nabla \sigma.
 \end{array}
 \]
 
 Do integration by parts, using boundary conditions $\eqref{e2.3}_{5,6,7}$, we have 
 \begin{eqnarray}
   \frac{1}{k} \int_{\Omega}   c^m (c^m-c^{m-1}) +\alpha \Vert  \nabla c^m   \Vert _{L^2}^2 +\frac{\alpha}{b}
 \Vert  \nabla_{\tau} c^m   \Vert _{L^2}^2 +\frac{\alpha}{bk} \int_{\partial\Omega}   c_{\tau} ^m 
 (c_{\tau} ^m-c_{\tau} ^{m-1})\nonumber \\
 =  - \int_{\Omega} c^m n^m f(c^m),\label{e4.1}\\
  \frac{1}{k} \int_{\Omega}   n^m (n^m-n^{m-1})+\beta \Vert  \nabla n^m   \Vert _{L^2}^2 -
 \int_{\Omega} g(n^m,\ c^m)  \nabla c^m\nabla n^m=0,\label{e4.2}\\
 \frac{1}{k} \int_{\Omega}   u^m (u^m-u^{m-1})+\xi \Vert  \nabla u^m   \Vert _{L^2}^2=\int_{\Omega} 
 u^m n^m \nabla \sigma. \label{e4.3}
\end{eqnarray}

Using relation (2.1) in the equations \eqref{e4.1} and \eqref{e4.2}, we get,
 \begin{equation}\label{e4.4}
 \begin{array}{ll}
\   &  \frac{1}{2k}( \Vert   c^m   \Vert _{L^2}^2  -  \Vert   c^{m-1}   \Vert _{L^2}^2  + \Vert  c^m-c^{m-1}  \Vert _{L^2}^2 )+
 \alpha \Vert  \nabla c^m   \Vert _{L^2}^2 +\frac{\alpha}{b} \Vert  \nabla_{\tau} c^m   \Vert _{L^2}^2 \\
 \  &  \hspace{0.5cm}   +\frac{\alpha}{2bk} 
 ( \Vert   c_{\tau}^m   \Vert _{L^2}^2  -  \Vert   c_{\tau}^{m-1}   \Vert _{L^2}^2  + \Vert  c_{\tau}^m-c_{\tau}^{m-1}  \Vert _{L^2}^2 )\\
  \   &  \hspace{1.5cm} \le \frac{f_1}{2} ( \Vert   c^m   \Vert _{L^2}^2 +  \Vert   n^m   \Vert _{L^2}^2 ).
 \end{array}
 \end{equation}
  
 \begin{equation}\label{e4.5}
 \begin{array}{ll}  
 \   &  \frac{1}{2k}( \Vert   n^m   \Vert _{L^2}^2  -  \Vert   n^{m-1}   \Vert _{L^2}^2  + \Vert  n^m-n^{m-1}  \Vert _{L^2}^2 )+
 \beta \Vert  \nabla n^m   \Vert _{L^2}^2\\
  \   & \hspace{1.5cm} \le g_1 ( \frac{1}{4\delta} \Vert  \nabla  c^m   \Vert _{L^2}^2 + \delta \Vert   \nabla n^m   \Vert _{L^2}^2 ).
\end{array}
\end{equation}

 Multiply \eqref{e4.4} by $\frac{g_1}{4\delta},\ $ \eqref{e4.5} by $\alpha$,  add the resulting equations, we get
 \begin{equation}\label{e4.6}
 \begin{array}{ll}
 \   &  \frac{g_1}{8k\delta}( \Vert   c^m   \Vert _{L^2}^2  -  \Vert   c^{m-1}   \Vert _{L^2}^2  + \Vert  c^m-c^{m-1}  \Vert _{L^2}^2 )+
 \frac{g_1\alpha}{4b\delta}
  \Vert  \nabla_{\tau} c^m   \Vert _{L^2}^2\\
 \   &  \hspace{0.5cm}  +\frac{g_1\alpha}{8kb\delta}  ( \Vert   c_{\tau}^m   \Vert _{L^2}^2  -  \Vert   c_{\tau}^{m-1}   \Vert _{L^2}^2  + \Vert  c_{\tau}^m-c_{\tau}^{m-1}  \Vert _{L^2}^2 )  \\
 \   &  \hspace{0.5cm}  + \frac{\alpha}{2k}( \Vert   n^m   \Vert _{L^2}^2  -  \Vert   n^{m-1}   \Vert _{L^2}^2  + \Vert  n^m-n^{m-1}  \Vert _{L^2}^2 ) + (\alpha \beta - g_1\alpha\delta)\Vert   \nabla n^m   \Vert _{L^2}^2\\
 \   &  \hspace{1.5cm} \le  \frac{g_1}{4\delta}  \frac{f_1}{2} ( \Vert   c^m   \Vert _{L^2}^2 +  \Vert   n^m   \Vert _{L^2}^2 ).
\end{array}
\end{equation}

Choose $\delta$ such that $\delta < \frac{\beta}{g_1}$, then multiply \eqref{e4.6} by $\frac{8k\delta}{g_1}$, we get
 \begin{equation}\label{e4.7}
 \begin{array}{ll}
 \   &  ( \Vert   c^m   \Vert _{L^2}^2  -  \Vert   c^{m-1}   \Vert _{L^2}^2  + \Vert  c^m-c^{m-1}  \Vert _{L^2}^2 )+
 \frac{2k\alpha}{b} \Vert  \nabla_{\tau} c^m   \Vert _{L^2}^2\\
 \   &  \hspace{0.5cm}  +\frac{\alpha}{b}  ( \Vert   c_{\tau}^m   \Vert _{L^2}^2  -  \Vert   c_{\tau}^{m-1}   \Vert _{L^2}^2  + \Vert  c_{\tau}^m-c_{\tau}^{m-1}  \Vert _{L^2}^2 )  \\
 \   &  \hspace{0.5cm}  + \frac{4\delta\alpha}{g_1}( \Vert   n^m   \Vert _{L^2}^2  -  \Vert   n^{m-1}   \Vert _{L^2}^2  + \Vert  n^m-n^{m-1}  \Vert _{L^2}^2 ) + \frac{8k\delta}{g_1}    (\alpha \beta - g_1\alpha\delta)\Vert   \nabla n^m   \Vert _{L^2}^2\\
 \   &  \hspace{1.5cm} \le  k  f_1 ( \Vert   c^m   \Vert _{L^2}^2 +  \Vert   n^m   \Vert _{L^2}^2 ).
\end{array}
\end{equation}

Let $\tilde\alpha := \min\{ \frac{\alpha}{b},\   \frac{4\delta\alpha}{g_1},\    \frac{8\delta}{g_1}    (\alpha \beta - g_1\alpha\delta)    \},\ $ 
then from \eqref{e4.7} , we have
\[ 
 \begin{array}{ll}
 \   &  ( \Vert   c^m   \Vert _{L^2}^2  -  \Vert   c^{m-1}   \Vert _{L^2}^2  + \Vert  c^m-c^{m-1}  \Vert _{L^2}^2 )+
k\tilde\alpha \Vert  \nabla_{\tau} c^m   \Vert _{L^2}^2\\
 \   &  \hspace{0.5cm}  +\tilde\alpha  ( \Vert   c_{\tau}^m   \Vert _{L^2}^2  -  \Vert   c_{\tau}^{m-1}   \Vert _{L^2}^2  + \Vert  c_{\tau}^m-c_{\tau}^{m-1}  \Vert _{L^2}^2 )  \\
 \   &  \hspace{0.5cm}  + \tilde\alpha( \Vert   n^m   \Vert _{L^2}^2  -  \Vert   n^{m-1}   \Vert _{L^2}^2  + \Vert  n^m-n^{m-1}  \Vert _{L^2}^2 ) + \tilde\alpha k\Vert   \nabla n^m   \Vert _{L^2}^2\\
 \   &  \hspace{1.5cm} \le  k  f_1 ( \Vert   c^m   \Vert _{L^2}^2 +  \Vert   n^m   \Vert _{L^2}^2 ).
\end{array}
\]

Summing the above inequality for $m=1,\ 2,\ 3,\  \cdots ,\ r,\quad  1\le r \le  N,\ $ we find
\[ 
 \begin{array}{ll}
 \   &   \Vert   c^r   \Vert _{L^2}^2    + \sum_{m=1}^r  \Vert  c^m-c^{m-1}  \Vert _{L^2}^2 +
k \tilde{\alpha}  \sum_{m=1}^r  \Vert  \nabla_{\tau} c^m   \Vert _{L^2}^2\\
 \   &  \hspace{0.5cm}  + \tilde{\alpha}  ( \Vert   c_{\tau}^r   \Vert _{L^2}^2+ \sum_{m=1}^r  
 \Vert  c_{\tau}^m-c_{\tau}^{m-1}  \Vert _{L^2}^2 )  \\  
 \   &  \hspace{0.5cm}  +  \tilde{\alpha}( \Vert   n^r   \Vert _{L^2}^2  + \sum_{m=1}^r  \Vert  n^m-n^{m-1}  \Vert _{L^2}^2 ) 
  + \tilde{\alpha} k  \sum_{m=1}^r \Vert   \nabla n^m   \Vert _{L^2}^2\\
 \   &  \hspace{1.5cm} \le  k  f_1 \sum_{m=1}^r ( \Vert   c^m   \Vert _{L^2}^2 +  \Vert   n^m   \Vert _{L^2}^2 )
 +  \Vert   c^0  \Vert _{L^2}^2 + \tilde{\alpha} \Vert   c_{\tau}^0  \Vert _{L^2}^2 +  \tilde{\alpha}\Vert   n^0   \Vert _{L^2}^2.
\end{array}
\]

From here, we derive the following inequality:
\begin{equation}\label{e4.8}
 \begin{array}{ll}
 \   &   \Vert   c^r   \Vert _{L^2}^2    + \sum_{m=1}^r  \Vert  c^m-c^{m-1}  \Vert _{L^2}^2 +
k   \sum_{m=1}^r  \Vert  \nabla_{\tau} c^m   \Vert _{L^2}^2\\
 \   &  \hspace{0.5cm}  +   \Vert   c_{\tau}^r   \Vert _{L^2}^2+ \sum_{m=1}^r  
 \Vert  c_{\tau}^m-c_{\tau}^{m-1}  \Vert _{L^2}^2   \\  
 \   &  \hspace{0.5cm}  +  \Vert   n^r   \Vert _{L^2}^2  + \sum_{m=1}^r  \Vert  n^m-n^{m-1}  \Vert _{L^2}^2  
  + k  \sum_{m=1}^r \Vert   \nabla n^m   \Vert _{L^2}^2\\
 \   &  \hspace{0.5cm} \le M( \Vert   c^0  \Vert _{L^2}^2 +  \Vert   c_{\tau}^0  \Vert _{L^2}^2 +\Vert   n^0   \Vert _{L^2}^2)
 +Mk \sum_{m=1}^r ( \Vert   c^m   \Vert _{L^2}^2 +  \Vert   n^m   \Vert _{L^2}^2 ).
 \end{array}
\end{equation}

Using Lemma 2.1 (discrete Gronwall's Lemma), we have
\begin{equation}\label{e4.9}
 \begin{array}{ll}
\  & \max_{1 \le r \le N} \Vert   c^r  \Vert _{L^2}^2 +  \max_{1 \le r \le N} \Vert   n^r  \Vert _{L^2}^2\\
  \   &  \hspace{0.5cm}\le M ( \Vert   c^0  \Vert _{L^2}^2 +  \Vert   c_{\tau}^0  \Vert _{L^2}^2 +\Vert   n^0   \Vert _{L^2}^2).
 \end{array}
\end{equation}

\[
\begin{array}{ll}
 \   &   k   \sum_{m=1}^r ( \Vert  \nabla_{\tau} c^m   \Vert _{L^2}^2 + \Vert   \nabla n^m   \Vert _{L^2}^2) \\
 \   &  \hspace{0.5cm}\le M ( \Vert   c^0  \Vert _{L^2}^2 +  \Vert   c_{\tau}^0  \Vert _{L^2}^2 +\Vert   n^0   \Vert _{L^2}^2)
 +Mkr\left( M ( \Vert   c^0  \Vert _{L^2}^2 +  \Vert   c_{\tau}^0  \Vert _{L^2}^2 +\Vert   n^0   \Vert _{L^2}^2)    \right),\\
  \   &  \hspace{0.5cm}\le (M+M^2 T)( \Vert   c^0  \Vert _{L^2}^2 +  \Vert   c_{\tau}^0  \Vert _{L^2}^2 +\Vert   n^0   \Vert _{L^2}^2).
\end{array}
\]

Write $M+M^2 T$ again as $M$, we have
\begin{equation}\label{e4.10}
\begin{array}{ll}
\   &   k   \sum_{m=1}^r ( \Vert  \nabla_{\tau} c^m   \Vert _{L^2}^2 + \Vert   \nabla n^m   \Vert _{L^2}^2) \\
 \   &  \hspace{0.5cm}\le M( \Vert   c^0  \Vert _{L^2}^2 +  \Vert   c_{\tau}^0  \Vert _{L^2}^2 +\Vert   n^0   \Vert _{L^2}^2).
\end{array}
\end{equation}

Similarly, we have
\begin{equation}\label{e4.11}
\begin{array}{ll}
\   &   \sum_{m=1}^r  (\Vert  c^m-c^{m-1}  \Vert _{L^2}^2  + \Vert  c_{\tau}^m-c_{\tau}^{m-1}  \Vert _{L^2}^2  +
 \Vert  n^m-n^{m-1}  \Vert _{L^2}^2  )\\
 \   &  \hspace{0.5cm}\le M( \Vert   c^0  \Vert _{L^2}^2 +  \Vert   c_{\tau}^0  \Vert _{L^2}^2 +\Vert   n^0   \Vert _{L^2}^2).
\end{array}
\end{equation}

\begin{equation}\label{e4.12}
 \max_{1 \le r \le N} \Vert   c_{\tau}^r  \Vert _{L^2}^2  \le M( \Vert   c^0  \Vert _{L^2}^2 +  \Vert   c_{\tau}^0  \Vert _{L^2}^2 +\Vert   n^0   \Vert _{L^2}^2).
 \end{equation}

 \begin{equation}\label{e4.13}
 \sum_{m=1}^r  \Vert  c_{\tau}^m - c_{\tau}^{m-1}  \Vert _{L^2}^2  \le M( \Vert   c^0  \Vert _{L^2}^2 +  \Vert   c_{\tau}^0  \Vert _{L^2}^2 +\Vert   n^0   \Vert _{L^2}^2).
\end{equation}

For the estimate of $u^m$, we proceed in an analogous way, and obtain the following inequality from \eqref{e4.3}.
\[
\begin{array}{ll}
\   &  \Vert   u^m  \Vert _{L^2}^2 - \Vert   u^{m-1}  \Vert _{L^2}^2 + \Vert   u^m -  u^{m-1} \Vert _{L^2}^2
+ 2k\xi \Vert  \nabla u^m  \Vert _{L^2}^2\\
 \   &  \hspace{0.5cm}\le kC(  \Vert   u^m  \Vert _{L^2}^2 +  \Vert   n^m  \Vert _{L^2}^2  ).
\end{array}
\]

Summing the above inequality for $m=1,\ 2,\ 3,\  \cdots ,\ r,\quad  1\le r \le  N,\ $ we find
\[
\begin{array}{ll}
\   &   \Vert   u^r   \Vert _{L^2}^2 + \sum_{m=1}^r  \Vert  u^m-u^{m-1}  \Vert _{L^2}^2  + 2k\xi  \sum_{m=1}^r
\Vert  \nabla u^m  \Vert _{L^2}^2\\
 \   &  \hspace{0.5cm}\le kC \sum_{m=1}^r  (  \Vert   u^m  \Vert _{L^2}^2 +  \Vert   n^m  \Vert _{L^2}^2  )
 +  \Vert   u^0  \Vert _{L^2}^2\\
 \   &  \hspace{0.5cm}\le k C  \sum_{m=1}^r   \Vert   u^m  \Vert _{L^2}^2
 + krC \max _ {1\le m \le r} \Vert   n^m  \Vert _{L^2}^2 +   \Vert   u^0  \Vert _{L^2}^2\\
  \   &  \hspace{0.5cm}\le  M(\Vert   u^0  \Vert _{L^2}^2+ \Vert   c^0  \Vert _{L^2}^2 +  \Vert   c_{\tau}^0  \Vert _{L^2}^2 +\Vert   n^0   \Vert _{L^2}^2) + Mk \sum_{m=1}^r  \Vert  u^m \Vert _{L^2}^2 .
\end{array}
\]

By Lemma 2.1 (discrete Gronwall's Lemma), we have 
 \begin{equation}\label{e4.14}
 \max_{1 \le r \le N}  \Vert   u^r  \Vert _{L^2}^2  \le M( \Vert   c^0  \Vert _{L^2}^2 +  \Vert   c_{\tau}^0  \Vert _{L^2}^2 +\Vert   n^0   \Vert _{L^2}^2 + \Vert   u^0  \Vert _{L^2}^2).
 \end{equation}
 
 \begin{equation}\label{e4.15}
  \sum_{m=1}^r  \Vert  u^m - u^{m-1}  \Vert _{L^2}^2  \le M( \Vert   c^0  \Vert _{L^2}^2 +  \Vert   c_{\tau}^0  \Vert _{L^2}^2 +\Vert   n^0   \Vert _{L^2}^2 + \Vert   u^0  \Vert _{L^2}^2).
 \end{equation}

 \begin{equation}\label{e4.16}
  k\sum_{m=1}^r \Vert \nabla u^m \Vert _{L^2}^2  \le M( \Vert   c^0  \Vert _{L^2}^2 +  \Vert   c_{\tau}^0  \Vert _{L^2}^2 +\Vert   n^0   \Vert _{L^2}^2 + \Vert   u^0  \Vert _{L^2}^2).
\end{equation}

 From \eqref{e4.9} - \eqref{e4.16}, we obtain the following estimates for Rothe functions and step functions:
  \begin{equation}\label{e4.17}
 \begin{array}{ll}
 \   &  \Vert (c_k,\ n_k) \Vert_{\left(  L^2(0,\ T;\ L^2(\Omega)) \right)^2} \le M; \hspace{0.5cm} \Vert (c_k,\ n_k) \Vert_{\left(  L^2(0,\ T;\ H^1(\Omega)) \right)^2} \le M ,\\
  \   &  \Vert u_k \Vert_{L^2(0,\ T;\ H)} \le M; \hspace{0.5cm}  \Vert u_k \Vert_{L^2(0,\ T;\ V)} \le M ,\\
  \   &  \Vert (\tilde{c}_k,\  \tilde{n}_k) \Vert_{\left(  L^2(0,\ T;\ L^2(\Omega)) \right)^2} \le M; \hspace{0.5cm} \Vert  (\tilde{c}_k,\  \tilde{n}_k) \Vert_{\left(  L^2(0,\ T;\ H^1(\Omega)) \right)^2} \le M ,\\
 \   &  \Vert \tilde{u}_k \Vert_{L^2(0,\ T;\ H)} \le M; \hspace{0.5cm}  \Vert  \tilde{u}_k \Vert_{L^2(0,\ T;\ V)} \le M,\\
 \   &  \Vert c_{k\tau} \Vert_{L^2(0,\ T;\ L^2(\partial\Omega)) } \le M; \hspace{0.5cm} \Vert c_{k\tau} \Vert_{L^2(0,\ T;\ H^1(\partial\Omega)) } \le M,\\ 
 \   &   \Vert  \tilde{c}_{k\tau} \Vert_{L^2(0,\ T;\ L^2(\partial\Omega)) } \le M; \hspace{0.5cm} \Vert  \tilde{c}_{k\tau} \Vert_{L^2(0,\ T; \   H^1(\partial\Omega)) } \le M,\\
  \   &  \Vert (\tilde{c}_k - c_k,\ \tilde{n}_k - n_k) \Vert_{\left(  L^2(0,\ T;\ L^2(\Omega)) \right)^2} \le M k;  \hspace{0.5cm}
   \Vert \tilde{u}_k - u_k \Vert_{ L^2(0,\ T;\  H) } \le M k.
 \end{array}
\end{equation}

From the above estimates of Rothe functions and step functions, we conclude that, there exist subsequences, still denoted as
 $(c_k,\ c_{k\tau},\ n_k,\ u_k),\  (\tilde{c}_k,\  \tilde{c}_{k\tau},\ \tilde{n}_k,\ \tilde{u}_k)$  such that as $k\rightarrow 0,\ $, we have
  \begin{equation}\label{e4.18}
\begin{array}{ll}
\  & (c_k,\ c_{k\tau},\  n_k) \rightharpoonup (c,\ c_{\tau},\  n);\   (\tilde{c}_k,\  \tilde{c}_{k\tau},\  \tilde{n}_k) \rightharpoonup (\tilde{c},\  \tilde{c}_{\tau},\  \tilde{n}) \    \mbox{weakly in} \    L^2 (0,\ T;\ L^2(\Omega)),\\
\  & (c_k,\ c_{k\tau},\  n_k) \rightharpoonup (c,\ c_{\tau},\  n);\   (\tilde{c}_k,\  \tilde{c}_{k\tau},\  \tilde{n}_k) \rightharpoonup (\tilde{c},\  \tilde{c}_{\tau},\  \tilde{n}) \    \mbox{weakly in} \    L^2 (0,\ T;\ H^1(\Omega)),\\
\  & u_k \rightharpoonup u;\   \tilde{u}_k \rightharpoonup \tilde{u}, \    \mbox{weakly in} \    L^2 (0,\ T;\ H),\\
\  & u_k \rightharpoonup u;\   \tilde{u}_k \rightharpoonup \tilde{u}, \    \mbox{weakly in} \    L^2 (0,\ T;\ V),\\
\  & (c_k,\ c_{k\tau},\ n_k,\ u_k) = (\tilde{c}_k,\  \tilde{c}_{k\tau},\ \tilde{n}_k,\ \tilde{u}_k)\     \mbox{almost everywhere.}
\end{array}
\end{equation}
 Next, we want to show there exist subsequences of $(c_k,\ n_k,\ u_k)$, still denoted as $(c_k,\ n_k,\ u_k)$, such that as 
 $k\rightarrow 0,\ $ we have 
  \begin{equation}\label{e4.19}
 \begin{array}{ll}
\  & (c_k,\  n_k) \rightarrow (c,\  n)\       \mbox{strongly in} \    L^2 (0,\ T;\ L^2(\Omega)),\\
\  & u_k \rightarrow u\     \mbox{strongly in} \    L^2 (0,\ T;\ H).
\end{array}
\end{equation}
 To this end, we apply Lemma 2.2, set $Y=\left( H^1(\Omega)   \right)^2\times V,\    X= \left( L^2(\Omega)   \right)^2\times H,\ $
 and set $p=2.\ $  The embedding of $Y$ into $X$ is compact.  Let $G$ be the family of functions 
 $\  (\tilde{c}_k,\ \tilde{n}_k,\ \tilde{u}_k)$. Then $G$ is bounded in $L^2 (0,\ T;\ Y)\ $ and $L^2 (0,\ T;\ X)\ $ due to \eqref{e4.17}.  Assume $(c^0,\ n^0,\ u^0) \in \left(H^1  \right)^2\times V,\  $ we want to show (2.2) in Lemma 2.2 holds.
 
 For this, we rewrite \eqref{e2.3} using Rothe functions $(\tilde{c}_k,\ \tilde{c}_{k\tau},\  \tilde{n}_k,\ \tilde{u}_k)\ $ and step functions
$ (c_k,\  c_{k\tau},\  n_k,\ u_k)$, then multiply $\eqref{e2.3}_1$ by $\phi_1$, $\eqref{e2.3}_2$ by $\phi_2$, $\eqref{e2.3}_3$ by $\phi_3$,  where  $\phi_1, \   \phi_2$ are any
functions in $H^1(\Omega)$, and $\phi_3$ is any function in $V,$ integrate by parts, using boundary conditions, we obtain
 \begin{equation}\label{e4.20}
\begin{array}{ll}
\  &\int_{\Omega} \partial_t \tilde{c}_k(t)\phi_1 dx + \alpha \int_{\Omega} \nabla c_k(t) \nabla\phi_1 dx + \frac{\alpha}{b} \int_{\partial\Omega}  \partial_t \tilde{c}_{k\tau}(t)\phi_{1\tau} d\sigma +\\ 
\  &\quad   \frac{\alpha}{b}  \int_{\partial\Omega} \nabla_{\tau} c_k(t) \nabla_{\tau}\phi_1 d\sigma  +         \int_{\Omega} u_k(t)\nabla c_k(t) \phi_1 dx= \int_{\Omega} -n_k(t) f(c_k(t))\phi_1 dx,\\
\ & \int_{\Omega} \partial_t \tilde{n}_k(t)\phi_2 dx +  \int_{\Omega}     \left( \beta\nabla n_k(t) -  g ( n_k(t),\ c_k(t)  )   \nabla c_k(t)   \right)  \nabla\phi_2 dx \\
  \    &\quad  + \int_{\Omega} u_k(t)\nabla n_k(t) \phi_2 dx=0,\\    
\   & \int_{\Omega} \partial_t \tilde{u}_k(t)\phi_3 dx +   \int_{\Omega} \xi  \nabla u_k(t) \nabla\phi_3 dx +  \int_{\Omega} (u_k(t)\cdot\nabla) u_k(t) \phi_3 dx=\\
 \   &\quad \int_{\Omega} n_k(t) \nabla \sigma \phi_3 dx.
 \end{array}
\end{equation}

 We then integrate \eqref{e4.20} between $t$ and $t+a,\   t\in (0,\  T), \  a>0,$ we get
 \begin{equation}\label{e4.21}
\begin{array}{ll}
\ & \int_{\Omega} (\tilde{c}_k(t+a) - \tilde{c}_k(t)  ) \phi_1 dx +  \frac{\alpha}{b} \int_{\partial\Omega}  (\tilde{c}_{k\tau}(t+a) - \tilde{c}_{k\tau}(t))\phi_{1\tau} d\sigma      = -  \alpha \int_t^{t+a}    \int_{\Omega} \nabla c_k(s) \nabla\phi_1 dx ds\\
\  &\quad  -    \frac{\alpha}{b}  \int_t^{t+a}  \int_{\partial\Omega} \nabla_{\tau} c_k(s) \nabla_{\tau}\phi_1 d\sigma ds   -\int_t^{t+a}   \int_{\Omega} u_k(s)\nabla c_k(s) \phi_1 dx ds
  - \int_t^{t+a} \int_{\Omega} n_k(s) f(c_k(s))\phi_1 dx ds,\\ 
 \ & \int_{\Omega} (\tilde{n}_k(t+a) - \tilde{n}_k(t)  ) \phi_2 dx=   -  \int_t^{t+a}    \int_{\Omega}     \left( \beta\nabla n_k(s) - g (n_k(s),\ c_k(s)) \nabla c_k(s)   \right)  \nabla\phi_2 dx ds \\   
  \    &\quad - \int_t^{t+a}  \int_{\Omega} u_k(s)\nabla n_k(s) \phi_2 dx ds,\\    
 \  &\int_{\Omega} (\tilde{u}_k(t+a) - \tilde{u}_k(t)  ) \phi_3 dx =  -  \int_t^{t+a}  \int_{\Omega} \xi  \nabla u_k(s) \nabla\phi_3 dx ds
  - \int_t^{t+a} \int_{\Omega} (u_k(s)\cdot\nabla) u_k(s) \phi_3 dx ds\\
 \   &\quad + \int_t^{t+a} \int_{\Omega} n_k(s) \nabla \sigma \phi_3 dx ds.
 \end{array}
\end{equation}

In the above equations, let
 \[
\begin{array}{ll}
\ & \phi_1 = \tilde{c}_k(t+a) - \tilde{c}_k(t),\\
\ & \phi_{1\tau} = \tilde{c}_{k\tau}(t+a) - \tilde{c}_{k\tau}(t),\\
\ & \phi_2 = \tilde{n}_k(t+a) - \tilde{n}_k(t),\\
\ & \phi_3 = \tilde{u}_k(t+a) - \tilde{u}_k(t).
\end{array}
\]

Then integrate from $0$ to $T-a$, we get from $\eqref{e4.21}_1$,
\[
\begin{array}{ll}
\ & \int_0^{T-a} \Vert \tilde{c}_k(t+a) - \tilde{c}_k(t) \Vert_{L^2}^2  +   \frac{\alpha}{b}\int_0^{T-a} \Vert \tilde{c}_{k\tau}(t+a) - \tilde{c}_{k\tau}(t) \Vert_{L^2}^2= I_1 + I_2 + I_3+I_4,\\
\end{array}
\]

where 
\[
\begin{array}{ll}
\ &   I_1=  -  \alpha\int_0^{T-a} \int_t^{t+a}    \int_{\Omega} \nabla c_k(s) \nabla\phi_1 dx ds dt, \\
\ &   I_2=   - \frac{\alpha}{b}  \int_0^{T-a} \int_t^{t+a}  \int_{\partial\Omega} \nabla_{\tau} c_k(s) \nabla_{\tau}\phi_1 d\sigma ds dt,\\
\ &   I_3=  -  \int_0^{T-a} \int_t^{t+a}    \int_{\Omega} u_k(s)\nabla c_k(s) \phi_1 dx ds dt, \\
\ &   I_4=  -  \int_0^{T-a} \int_t^{t+a}    \int_{\Omega} n_k(s) f(c_k(s))\phi_1  dx ds dt. 
\end{array}
\]

Using H$\ddot{o}$lder's inequality, Fubini's Theorem, and \eqref{e4.17}, we have 
\[
\begin{array}{ll}
\vert I_1 \vert & \le \alpha\int_0^{T-a} \int_t^{t+a} \Vert \nabla c_k(s) \Vert _{L^2} \Vert \nabla  (\tilde{c}_k(t+a) - \tilde{c}_k(t)) \Vert _{L^2}
ds dt\\
\   & =    \alpha\int_0^{T-a} \Vert \nabla  (\tilde{c}_k(t+a) - \tilde{c}_k(t)) \Vert _{L^2} \left( \int_t^{t+a} \Vert \nabla c_k(s) \Vert _{L^2}  ds  \right) dt\\
\   & \le   a^{\frac{1}{2}}  \left(\int_0^{T} \Vert \nabla c_k(s) \Vert _{L^2}^2  ds\right)^ \frac{1}{2}    \alpha\int_0^{T-a} \Vert \nabla  (\tilde{c}_k(t+a) - \tilde{c}_k(t)) \Vert _{L^2} dt\\
\   & \le a^{\frac{1}{2}}  M    \alpha   {(T-a)}^{\frac{1}{2}} 2  \left( \int_0^{T} \Vert \nabla  \tilde{c}_k(t) \Vert _{L^2}^2 dt \right)^{\frac{1}{2}} \\
\   & \le a^{\frac{1}{2}} M    \alpha   {T}^{\frac{1}{2}}.
\end{array}
\]

Then we have $I_1 \rightarrow 0$ as $a \rightarrow 0.$ Similarly, we have $I_2 \rightarrow 0$ as $a \rightarrow 0.$

For the estimates of  $I_3,\ $by H$\ddot{o}$lder's inequality, Fubini's Theorem, Sobolev embedding theorem, and \eqref{e4.17}, we have
\[
\begin{array}{ll}
\vert I_3 \vert & \le \int_0^{T-a}\int_t^{t+a}\Vert \ u_k(s) \Vert_{L^6} \Vert \nabla c_k(s) \Vert _{L^2}  
 \Vert  \tilde{c}_k(t+a) - \tilde{c}_k(t)  \Vert _{L^3} dsdt\\
\  & \le \int_0^{T-a} 2  \Vert \nabla \tilde{c}_k(t)  \Vert _{L^2} \left\{ \int_t^{t+a}\Vert  u_k(s) \Vert_V \Vert \nabla c_k(s) \Vert _{L^2}
ds\right\}dt\\ 
\  & \le \int_0^{T} \Vert u_k(s) \Vert _V  \Vert \nabla c_k(s) \Vert _{L^2}   \left(   \int_{[s-a,\ s]\bigcap [0,\ T-a] }  
 2\Vert \nabla  \tilde{c}_k(t)   \Vert _{L^2} dt \right)  ds\\
\  & \le \int_0^{T} \Vert \ u_k(s) \Vert _V  \Vert \nabla c_k(s) \Vert _{L^2}  \left\{ 2a^{\frac{1}{2}}   \left(   \int_0^T 
 \Vert \nabla  \tilde{c}_k(t)   \Vert _{L^2}^2 dt   \right)^{\frac{1}{2}}      \right\} ds\\
\  & \le  a^{\frac{1}{2}} M   \int_0^{T} \Vert \ u_k(s) \Vert _V  \Vert \nabla c_k(s) \Vert _{L^2} ds\\
\  & \le  a^{\frac{1}{2}} M \left(  \int_0^{T} \Vert \ u_k(s) \Vert _V^2 ds  \right)^{\frac{1}{2}}     \left(  \int_0^{T}
   \Vert \nabla c_k(s) \Vert _{L^2}^2 ds\right)^ { \frac{1}{2} }\\
\  &\le a^{\frac{1}{2}} M.
\end{array}
\]

Then we have $I_3 \rightarrow 0$ as $a \rightarrow 0.$

We now continue to estimate $\vert I_4 \vert:$
\[
\begin{array}{ll}
\vert I_4 \vert & \le \int_0^{T-a} \int_t^{t+a} f_1  \Vert  n_k(s) \Vert _{L^2} \Vert   \tilde{c}_k(t+a) - \tilde{c}_k(t) \Vert _{L^2} ds dt\\
\  &= f_1  \int_0^{T-a}   \Vert   \tilde{c}_k(t+a) - \tilde{c}_k(t) \Vert _{L^2}  \left(\int_t^{t+a}\Vert  n_k(s) \Vert _{L^2} ds\right) dt\\
\   & \le f_1 a^{\frac{1}{2}}  \left(\int_0^{T}\Vert  n_k(s) \Vert _{L^2}^2 ds\right)^{\frac{1}{2}}    \int_0^{T-a}   \Vert   \tilde{c}_k(t+a) - \tilde{c}_k(t) \Vert _{L^2} dt\\
\   & \le f_1 a^{\frac{1}{2}}  {(T-a)}^{\frac{1}{2}} M.
\end{array}
\]

We have $I_4 \rightarrow 0$ as $a \rightarrow 0.$

Therefore, we have proved that
 \begin{equation}\label{e4.22}
\begin{array}{ll}
\ &\int_0^{T-a} \Vert \tilde{c}_k(t+a) - \tilde{c}_k(t) \Vert_{L^2}^2 \  \rightarrow 0  \quad \mbox{as}\     a \rightarrow 0.
\end{array}
\end{equation}

Proceeding in an analogous way, we have
\begin{eqnarray}
&\  &\int_0^{T-a} \Vert \tilde{n}_k(t+a) - \tilde{n}_k(t) \Vert_{L^2}^2 \  \rightarrow 0  \quad \mbox{as}\     a \rightarrow 0. \label{e4.23}\\ 
&\  &\int_0^{T-a} \Vert \tilde{u}_k(t+a) - \tilde{u}_k(t) \Vert_{L^2}^2 \  \rightarrow 0   \quad \mbox{as}\     a \rightarrow 0.\label{e4.24} 
\end{eqnarray}

With \eqref{e4.22} - \eqref{e4.24}, and lemma 2.2, we conclude that \eqref{e4.19} holds.

Proceed in the same way as we derive \eqref{e3.28} - \eqref{e3.29}, we also have,
\begin{eqnarray}
&\  & f(c_k) \rightarrow f(c) \quad  \mbox{strongly in} \   L^2(0,\ T,\  L^2(\Omega)),\label{e4.26}\\
&\  & g(n_k,\  c_k) \rightarrow g(n,\  c) \quad  \mbox{strongly in} \   L^2(0,\ T,\  L^2(\Omega)).\label{e4.27}
\end{eqnarray}

Next, we want to pass the limit in \eqref{e4.20} by letting $k\rightarrow 0,$ to this end, we consider
$\phi_1,\   \phi_2\in  C^{\infty}\bigcap   H^1,\   \phi_3\in  \Upsilon,$ and $(\Psi_1,\  \Psi_2,\  \Psi_3 ) \in {\left(   C^1([0,\ T],\ R)   \right)}^3,$ with $\Psi_1(T)=\Psi_2(T)=\Psi_3(T)=0.$

Multiply $ \eqref{e4.20}_1$ by $\Psi_1$, $ \eqref{e4.20}_2$ by $\Psi_2$, $ \eqref{e4.20}_3$ by $\Psi_3$, integrate over $[0,\ T],\ $ integration by parts, we get
\begin{equation}\label{e4.28}
\begin{array}{ll}
\  & -\int_0^T       \int_{\Omega}  \tilde{c}_k(t)\phi_1  \Psi_1'(t)  dx dt  + \alpha\int_0^T     \int_{\Omega} \nabla c_k(t) \nabla\phi_1  \Psi_1(t) 
 dx dt -\frac{\alpha}{b} \int_0^T    \int_{\partial\Omega}  \tilde{c}_{k\tau}(t)\phi_{1\tau}  \Psi_1'(t)  d\sigma dt \\
 \  &\quad +\frac{\alpha}{b} \int_0^T     \int_{\partial\Omega} \nabla_{\tau} c_k(t) \nabla_{\tau}\phi_1  \Psi_1(t) d\sigma dt 
 + \int_0^T  \int_{\Omega} u_k(t)\nabla c_k(t) \phi_1 \Psi_1(t) dx  dt
= \Psi_1(0) \int_{\Omega}  \tilde{c}_k(x,\ 0)\phi_1   dx \\
 \  &\quad  + \frac{\alpha}{b} \Psi_1(0) \int_{\partial\Omega}  \tilde{c}_{k\tau}(x,\ 0)\phi_{1\tau}   d\sigma
 -  \int_0^T   \int_{\Omega} n_k(t) f(c_k(t))\phi_1 \Psi_1(t) dx dt,\\
\  & -\int_0^T       \int_{\Omega}  \tilde{n}_k(t)\phi_2  \Psi_2'(t)  dx dt   + \int_0^T  \int_{\Omega}     \left( \beta\nabla n_k(t) - g (n_k(t),\ c_k(t)) \nabla c_k(t)   \right)  \nabla\phi_2 \Psi_2(t) dx dt\\
\  & \hspace{1cm} + \int_0^T  \int_{\Omega} u_k(t)\nabla n_k(t) \phi_2  \Psi_2(t)  dx dt =   \Psi_2(0)  \int_{\Omega}  \tilde{n}_k(x,\ 0)\phi_2   dx,\\    
 \  & -\int_0^T     \int_{\Omega} \tilde{u}_k(t)\phi_3   \Psi_3'(t)  dx dt   + \int_0^T  \int_{\Omega} \xi  \nabla u_k(t) \nabla\phi_3 
  \Psi_3(t) dx dt + \int_0^T  \int_{\Omega} (u_k(t)\cdot\nabla) u_k(t) \phi_3 \Psi_3(t) dx dt\\
  \ &\hspace{1cm}=   \Psi_3(0)  \int_{\Omega}  \tilde{u}_k(x,\ 0)\phi_3 dx   + \int_0^T \int_{\Omega} n_k(t) \nabla \sigma \phi_3 \Psi_3(t) dx dt.
 \end{array}
\end{equation}

Now we take the limit in \eqref{e4.28} by letting $k \rightarrow 0,\ $ and using \eqref{e4.18}, \eqref{e4.19}, \eqref{e4.26}, \eqref{e4.27}, we get
\begin{equation}\label{e4.29}
\begin{array}{ll}
\ & -\int_0^T       \int_{\Omega}  c(t)\phi_1  \Psi_1'(t)  dx dt  + \alpha\int_0^T     \int_{\Omega} \nabla c(t) \nabla\phi_1  \Psi_1(t) 
 dx dt  -\frac{\alpha}{b} \int_0^T    \int_{\partial\Omega}  c_{\tau}(t)\phi_{1\tau}  \Psi_1'(t)  d\sigma dt \\
\ & +\frac{\alpha}{b} \int_0^T  \int_{\partial\Omega} \nabla_{\tau}  c(t) \nabla_{\tau}\phi_1  \Psi_1(t) d\sigma dt 
 + \int_0^T  \int_{\Omega} u(t)\nabla c(t) \phi_1 \Psi_1(t) dx  dt\\
\ & = \Psi_1(0) \int_{\Omega}  c(x,\ 0)\phi_1   dx 
+ \frac{\alpha}{b} \Psi_1(0) \int_{\partial\Omega}  c_{\tau}(x,\ 0)\phi_{1\tau}   d\sigma
-  \int_0^T   \int_{\Omega} n(t) f(c(t))\phi_1 \Psi_1(t) dx dt,\\
\ & -\int_0^T       \int_{\Omega}  n(t)\phi_2  \Psi_2'(t)  dx dt   + \int_0^T  \int_{\Omega}     \left( \beta\nabla n(t) - g (n(t),\ c(t)) \nabla c(t)   \right)  \nabla\phi_2 \Psi_2(t) dx dt\\
 \ & + \int_0^T  \int_{\Omega} u(t)\nabla n(t) \phi_2  \Psi_2(t)  dx dt =   \Psi_2(0)  \int_{\Omega}  n(x,\ 0)\phi_2   dx,\\    
\ &-\int_0^T     \int_{\Omega} u(t)\phi_3   \Psi_3'(t)  dx dt   + \int_0^T  \int_{\Omega} \xi  \nabla u(t) \nabla\phi_3 
  \Psi_3(t) dx dt + \int_0^T  \int_{\Omega} (u(t)\cdot\nabla) u(t) \phi_3 \Psi_3(t) dx dt\\
\ &=   \Psi_3(0)  \int_{\Omega}  u(x,\ 0)\phi_3 dx   + \int_0^T \int_{\Omega} n(t) \nabla \sigma\phi_3 \Psi_3(t) dx dt.
\end{array}
\end{equation}

\eqref{e4.29} holds for any $\phi_1,\   \phi_2\in  C^{\infty}\bigcap   H^1,\   \phi_3\in  \Upsilon,$ by continuity, it holds for any
 $\phi_1,\   \phi_2\in H^1,\   \phi_3\in  V$.

By choosing $(\Psi_1,\  \Psi_2,\  \Psi_3 ) \in {\left(   C_0^{\infty}[0,\ T]   \right)}^3,\ $ we draw conclusion that (1.6) holds, in 
the weak sense on $(0,\ T).$ By standard argument, we have $ (c(0),\ n(0),\ u(0)) = (c_0,\  n_0,\  u_0).$ The proof of Theorem 1.1
is completed.

\end{document}